# Nonsymmetric Koornwinder polynomials and duality

By Siddhartha Sahi*

## 1. Introduction

In the fundamental work of Lusztig [L] on affine Hecke algebras, a special role is played by the root system of type $\widetilde{C}_n$. The affine Hecke algebra is a deformation of the group algebra of an affine Weyl group which usually depends on as many parameters as there are distinct root lengths, i.e. one or two for an irreducible root system. However in the $\widetilde{C}_n$ case, the Hecke algebra $H$ has three parameters, corresponding to the fact that there is a simple coroot which is divisible by 2.

Recently, Cherednik [C1]–[C3] has introduced the notion of a double affine Hecke algebra, and has used it to prove several conjectures on Macdonald polynomials. These polynomials, and Cherednik's double affine Hecke algebra, involve two or three parameters, i.e., one more than the number of root lengths.

In this paper, motivated by the work of Noumi [N], we define a double affine Hecke algebra for the $\widetilde{C}_n$ case, which depends on three additional parameters, making six altogether. The associated orthogonal polynomials are precisely the six-parameter family of polynomials $P_\lambda$ introduced by Koornwinder in [Ko].

These polynomials are themselves quite remarkable. *Every* symmetric Macdonald polynomial [M] associated to a classical root system (i.e. those of types $A$, $B$, $C$, $D$, and the two classes of type $BC$, not considered by Cherednik) can be obtained from the $P_\lambda$ by a suitable limiting procedure [D]. Moreover, for $n = 1$, the $P_\lambda$ become the Askey-Wilson polynomials, which sit atop an impressive hierarchy of orthogonal polynomials in one variable [AW].

Koornwinder and Macdonald have formulated several conjectures for these polynomials, which are analogous to those proved by Cherednik for Macdonald polynomials. These are the "constant term," "norm," "evaluation," and

*This work was supported by an NSF grant, and carried out in part during the author's visit to Japan. The author wishes to thank M. Wakayama, H. Ochiai, K. Mimachi for their hospitality, and M. Noumi for explaining the results of [N].

"duality" conjectures. For a certain five-parameter subfamily of the $P_\lambda$, these conjectures were proved by van Diejen in [D].

In the general (six-parameter) setting, van Diejen has shown that either of the last two conjectures implies the other three. One of the results in this paper is a proof of the duality conjecture which implies all the rest by van Diejen's work.

Here is an outline of the paper: After a brief summary of the relevant results of Koornwinder and Noumi, we define the six-parameter double affine Hecke algebra $\mathcal{H}$ and establish its basic properties, including the existence of an involution. Next, we introduce certain commutators $S_i$ in $\mathcal{H}$, called the *intertwiners*, and use them to construct a family of polynomials $\{E_\alpha\}$. We call these the nonsymmetric Koornwinder polynomials, and we describe their relationship to the $P_\lambda$. Finally, we establish the duality conjecture for $P_\lambda$ together with its analog for $E_\alpha$.

A substantial part of this paper is directly motivated by the results of Cherednik in the two-/three-parameter setting. The idea of using intertwiners as creation operators was introduced in [K], [KS] and [S] for $\text{GL}_n$, and in [C4] for other root systems.

We have avoided one layer of notational complexity by identifying the coroot lattice of $C_n$ with $\mathbb{Z}^n$. Thus we have suppressed explicit reference to roots and weights. Implicitly, though, these are ubiquitous.

We have also obtained fairly precise results concerning the orthogonality and triangularity of the nonsymmetric Koornwinder polynomials, which we shall report elsewhere.

Finally, we remark that according to the note added in proof to [D], Macdonald has informed van Diejen that he has proved the evaluation conjecture. By van Diejen's work, this would also imply the duality conjecture for the $P_\lambda$, though not for the $E_\alpha$.

## 2. Preliminaries

In this section we briefly recall some results of Koornwinder, Lusztig, and Noumi which we shall need. For more details the reader should consult [Ko], [L], and [N].

We fix six indeterminates $q, t, t_0, t_n, u_0, u_n$, and let $\mathbb{F}$ be the field of rational functions in their *square roots*. We also define

(1) $\quad a = t_n^{1/2} u_n^{1/2}, b = -t_n^{1/2} u_n^{-1/2}, c = q^{1/2} t_0^{1/2} u_0^{1/2}, d = -q^{1/2} t_0^{1/2} u_0^{-1/2}.$

Let $\mathcal{R} = \mathbb{F}[x_1^{\pm 1}, \cdots, x_n^{\pm 1}]$ be the ring of Laurent polynomials in $n$ variables over the field $\mathbb{F}$, and let $\mathcal{S}$ be the subring consisting of *symmetric* polynomials, i.e. those which are invariant under permutations and inversions of the variables.

2.1. *Koornwinder polynomials.* In [Ko], Koornwinder defined a basis $\{P_\lambda\}$ of $\mathcal{S}$ which is indexed by $\lambda \in \mathbb{Z}^n$ with $\lambda_1 \geq \cdots \geq \lambda_n \geq 0$, and defined as follows:

Let $T_{q,x_i}$ denote the $i^{\text{th}}$ $q$-shift operator acting on $\mathcal{R} := \mathbb{F}[x_1^{\pm 1}, \cdots, x_n^{\pm 1}]$ by

$$T_{q,x_i} f(x_1, \cdots, x_i, \cdots, x_n) := f(x_1, \cdots, qx_i, \cdots, x_n).$$

Consider the following $q$-difference operator

$$D := \sum_{i=1}^n \Phi_i(x)(T_{q,x_i} - 1) + \sum_{i=1}^n \Phi_i(x^{-1})(T_{q,x_i}^{-1} - 1)$$

where $\Phi_i(x)$ is a rational function in $x_1, \cdots, x_n$ defined by

$$\Phi_i(x) := \frac{(1 - ax_i)(1 - bx_i)(1 - cx_i)(1 - dx_i)}{(1 - x_i^2)(1 - qx_i^2)} \prod_{\substack{j=1 \\ j \neq i}}^n \frac{\left(1 - tx_i x_j^{-1}\right)(1 - tx_i x_j)}{\left(1 - x_i x_j^{-1}\right)(1 - x_i x_j)}.$$

Koornwinder showed that $D$ preserves $\mathcal{S}$ and is diagonalizable with distinct eigenvalues

$$d_\lambda = \sum_{i=1}^n \left[q^{-1}abcd\, t^{2n-i-1}(q^{\lambda_i} - 1) + t^{i-1}(q^{-\lambda_i} - 1)\right].$$

The Koornwinder polynomial $P_\lambda$ is characterized by the equation

(2) $$DP_\lambda = d_\lambda P_\lambda,$$

together with the condition that the coefficient of $x^\lambda := x_1^{\lambda_1} \cdots x_n^{\lambda_n}$ in $P_\lambda$ is 1.

It turns out that Koornwinder's operator $D$ is one among a commuting family of difference operators, all of which are simultaneously diagonalized by the $P_\lambda$. These higher operators were constructed abstractly by Noumi, and explicitly by van Diejen.

To describe the results of Noumi, we need to introduce some additional notation.

2.2. *The affine Weyl group.* The affine Weyl group $W$ of type $\widetilde{C}_n$ is generated by elements $s_0, s_1, \cdots, s_n$ which satisfy $s_i^2 = 1$ and, for $n > 1$, also satisfy the braid relations

$$s_i s_j s_i \cdots = s_j s_i s_j \cdots$$

with two, three, or four terms on each side accordingly as $i$ and $j$ are connected by zero, one or two lines in the Coexeter graph

$$0 = 1 - 2 - \cdots\cdots\cdots - (n-2) - (n-1) = n.$$

The (finite) Weyl group $W_0$ of type $C_n$ is the subgroup generated by $s_1, \cdots, s_n$.

$W$ has a natural faithful affine action on $V = \mathbb{R}^n$ in which

(3) $$s_n \cdot v = (v_1, \cdots, v_{n-1}, -v_n), \quad s_0 \cdot v = (-v_1 - 1, v_2, \cdots, v_n)$$

and the other $s_i$ act by interchanging $v_i$ and $v_{i+1}$. Indeed $W$ is the semidirect product of $W_0$ and $\tau(\mathbb{Z}^n)$, where $\tau^v$ denotes the translation by $v$. In terms of the generators,

$$\tag{4} \tau^v = \tau_1^{v_1} \cdots \tau_n^{v_n}, \quad \tau_i = (s_i \cdots s_{n-1})(s_n \cdots s_0)(s_1^{-1} \cdots s_{i-1}^{-1}).$$

Since the $W$-action on $V$ is affine, we get a representation of $W$ on the space $\widetilde{V}$ of affine linear functionals on $V$, which we identify with $V \times \mathbb{R}\,\delta$ as follows:

$$\tag{5} \langle v + r\delta, v' \rangle := v_1' v_1 + \cdots + v_n' v_n + r; \quad \langle w(v + r\delta), v' \rangle := \langle v + r\delta, w^{-1} \cdot v' \rangle.$$

This representation is given explicitly by

$$\tag{6} \begin{aligned} s_i(v + r\delta) &= s_i v + r\delta, \quad i \neq 0 \\ s_0(v + r\delta) &= (-v_1, v_2, \cdots, v_n) + (r - v_1)\delta. \end{aligned}$$

We define an exponential map from the lattice $\mathbb{Z}^n \times \mathbb{Z}\delta \subset \widetilde{V}$ to $\mathcal{R}$ by

$$\tag{7} x^{v+k\delta} := q^{-k} x_1^{v_1} \cdots x_n^{v_n}, \quad v \in \mathbb{Z}^n, \ k \in \mathbb{Z}.$$

Since the lattice is preserved under (6), we get a representation of $W$ on $\mathcal{R}$ by putting

$$\tag{8} w(x^{\widetilde{v}}) := x^{w(\widetilde{v})}; \quad \widetilde{v} \in \mathbb{Z}^n \times \mathbb{Z}\delta.$$

Then $W$ acts by algebra homomorphisms and $\mathcal{S} = \mathcal{R}^{W_0}$. Explicitly we have

$$\tag{9} \begin{aligned} s_0 f(x) &= f(q x_1^{-1}, x_2, \cdots, x_n) \\ s_i f(x) &= f(x_1, \cdots, x_{i+1}, x_i, \cdots, x_n) \quad i \neq 0, n \\ s_n f(x) &= f(x_1, \cdots, x_{n-1}, x_n^{-1}) \\ \tau_i &= T_{q, x_i}. \end{aligned}$$

2.3. *The affine Hecke algebra.* The affine Hecke algebra $H$ of type $\widetilde{C}_n$ is generated over $\mathbb{F}$ by elements $T_0, T_1, \cdots, T_n$ which satisfy the same braid relations as the $s_i$, and also satisfy

$$T_i - T_i^{-1} = t_i^{1/2} - t_i^{-1/2}$$

where $t_1 = \cdots = t_{n-1} = t$.

The elements $T_1, \cdots, T_n$ generate the (finite) Hecke algebra $H_0$ of type $C_n$. $H$ and $H_0$ have natural bases $\{T_w\}$ consisting of $w$ in $W$ and $W_0$, respectively, where

$$\tag{10} T_w = T_{i_1} \cdots T_{i_l}$$

if $w = s_{i_1} \cdots s_{i_l}$ is a *reduced* (i.e. shortest) expression of $w$ in terms of the $s_i$.

The analogs of the translations $\tau_i$ in (4) are the elements

(11) $\qquad Y_i = (T_i \cdots T_{n-1})(T_n \cdots T_0)(T_1^{-1} \cdots T_{i-1}^{-1}), \quad i = 1, \cdots, n.$

Lusztig [L] has shown that the $Y_i$ commute pairwise and generate a subalgebra
$$\mathcal{R}_Y \approx \mathbb{F}[Y_1^{\pm 1}, \cdots, Y_n^{\pm 1}],$$
and that multiplication gives us a vector space isomorphism

(12) $\qquad\qquad\qquad H_0 \otimes \mathcal{R}_Y \approx H.$

2.4. *The Noumi representation.* Let $s_i$ act on $\mathcal{R}$ by (9); then Noumi in [N] has shown that the following map $\pi$ extends to a representation of $H$ on $\mathcal{R}$:

(13) $\qquad \pi(T_0^{\pm 1}) := t_0^{\pm 1/2} + t_0^{-1/2} \dfrac{(1 - cx_1^{-1})(1 - dx_1^{-1})}{1 - qx_1^{-2}}(s_0 - 1)$

$\qquad\qquad \pi(T_i^{\pm 1}) := t_i^{\pm 1/2} + t_i^{-1/2} \dfrac{(1 - t_i x_i x_{i+1}^{-1})}{(1 - x_i x_{i+1}^{-1})}(s_i - 1) \quad i \neq 0, n$

$\qquad\qquad \pi(T_n^{\pm 1}) := t_n^{\pm 1/2} + t_n^{-1/2} \dfrac{(1 - ax_n)(1 - bx_n)}{1 - x_n^2}(s_n - 1),$

where $a$, $b$, $c$, $d$ are as in (1).

Moreover if $f$ is in $\mathcal{S}$ then $\pi(f) := \pi(f(Y_1, \cdots, Y_n))$ preserves $\mathcal{S}$. Noumi showed that the restriction of $\pi(Y_1 + \cdots + Y_n + Y_1^{-1} + \cdots + Y_n^{-1})$ to $\mathcal{S}$ is a linear combination of the Koornwinder operator $D$ and a scalar. This means that the Koornwinder polynomials are simultaneous eigenfunctions for the $\pi(f)$. More precisely, $P_\lambda$ is characterized by

(14) $\qquad\qquad\qquad \pi(f)P_\lambda(x) = f(q^{\lambda + \rho})P_\lambda(x),$

where $q^{\lambda + \rho}$ means $(q^{\lambda_1 + \rho_1}, \cdots, q^{\lambda_n + \rho_n})$ and $\rho$ is defined by

(15) $\qquad q^{\rho_i} = st^{n-i}, \quad \text{with } s := (t_0 t_n)^{1/2} = \sqrt{q^{-1}abcd}.$

*Remark.* The fact that $\pi$ extends to a representation can be derived from Proposition 3.6 in [L], along the lines suggested in Proposition 4.6 of [M].

## 3. The double affine Hecke algebra

We now introduce the algebra $\mathcal{H}$ which will be the principal object of study in this paper. For convenience, we will write $Z \sim z$ as an abbreviation for the relation
$$Z - Z^{-1} = z^{1/2} - z^{-1/2}.$$

*Definition.* Let $\mathcal{H}$ be the algebra generated over $\mathbb{F}$ by elements $T_i^{\pm 1}$, $i = 0, \cdots n$, and commuting elements $X_i^{\pm 1}$, $i = 1, \cdots, n$, subject to the relations:

(i) $T_i \sim t_i$,
(ii) the $\widetilde{C}_n$ braid relations for the $T_i$'s,
(iii) $T_i X_j = X_j T_i$ if $|i - j| > 1$, or if $i = n$ and $j = n - 1$,
(iv) $T_i X_i = X_{i+1} T_i^{-1}$, $i = 1, \cdots, n-1$,
(v) $X_n^{-1} T_n^{-1} \sim u_n$,
(vi) $U_0 \equiv q^{-1/2} T_0^{-1} X_1 \sim u_0$.

If we set $u_0 = u_n = 1$ and $t_0 = t_n$, then $\mathcal{H}$ specializes to the three-parameter double affine Hecke algebra considered in [C1] for the affine root system $\widetilde{C}_n$.

Our definition is motivated by the following considerations:

Define a map $\pi$ from the generators of $\mathcal{H}$ to $\mathrm{End}(\mathcal{R})$ by letting $\pi(T_i^{\pm 1})$ be as in (13), and letting $\pi(X_i^{\pm 1})$ be the operator of multiplication by $x_i^{\pm 1}$.

3.1. THEOREM.   *The map $\pi$ extends to a representation of $\mathcal{H}$ on $\mathcal{R}$.*

*Proof.* We need only verify that (i)–(vi) hold for $\pi(T_i^{\pm 1})$ and $\pi(X_i^{\pm 1})$.

The relations (i) and (ii) follow from Noumi's result, and (iii) and (iv) are easily verified using the formulas. For (v) we have

$$X_n^{-1} T_n^{-1} \mapsto t_n^{-1/2} x_n^{-1} + t_n^{-1/2} \frac{(1 - ax_n)(1 - bx_n)}{1 - x_n^2}(x_n^{-1} s_n - x_n^{-1})$$

$$T_n X_n \mapsto t_n^{1/2} x_n + t_n^{-1/2} \frac{(1 - ax_n)(1 - bx_n)}{1 - x_n^2}(s_n x_n - x_n).$$

Since $s_n x_n = x_n^{-1} s_n$ and $x_n - x_n^{-1} = -x_n^{-1}(1 - x_n^2)$, we get

$$X_n^{-1} T_n^{-1} - T_n X_n \mapsto t_n^{-1/2} x_n^{-1} - t_n^{1/2} x_n - t_n^{-1/2} x_n^{-1}(1 - ax_n)(1 - bx_n)$$
$$= -(t_n^{1/2} + t_n^{-1/2} ab) x_n + t_n^{-1/2}(a + b).$$

Substituting $a = t_n^{1/2} u_n^{1/2}$ and $b = -t_n^{1/2} u_n^{-1/2}$, this becomes $u_n^{1/2} - u_n^{-1/2}$ proving (v).

Relation (vi) can be proved similarly using $s_0 x_1 = q x_1^{-1} s_0$.   □

3.2. THEOREM.   *The representation $\pi$ is faithful.*

*Proof.* We first note that in any word in $\mathcal{H}$ involving the generators, the relations (i)–(vi) allow us to commute the $T_i$'s past the $X_j$'s. Thus every element of $\mathcal{H}$ can be written as a linear combination of $X^\alpha T_w$, where $X_\alpha = X_1^{\alpha_1} \cdots X_n^{\alpha_n}$ and $T_w$ is as in (10).

Suppose a *nontrivial* linear combination maps to 0 under $\pi$. We then get
$$\sum c_{w,\alpha} x^\alpha \pi(T_w) = 0$$
in $\operatorname{End}(\mathcal{R})$, where $c_{w,\alpha}$ are scalars in $\mathbb{F}$, not all zero.

The left side is a rational expression in the square roots of $q, t, t_0, t_n, u_0, u_n$, and we consider what happens if we specialize the last five indeterminates to 1. By clearing denominators and eliminating common factors, we may ensure that at least *some* of the $c_{w,\alpha}$ have nonzero specializations; by (1), (10), and (13), $\pi(T_w)$ specializes to the action of $w$ as in (9).

Thus we get a nontrivial dependence relation in $\operatorname{End}(\mathcal{R})$ of the form
$$\sum c_{w,\alpha} x^\alpha w = 0, \quad w \in W, \alpha \in \mathbb{Z}^n.$$
Since $W = W_0 \tau(\mathbb{Z}^n)$, we can rewrite this as
$$\sum c_{w,\alpha,\beta} x^\alpha w \tau^\beta = 0, \quad w \in W_0, \alpha, \beta \in \mathbb{Z}^n.$$
Collecting the terms for $\beta$, we get
$$\sum x^\alpha w p_{w,\alpha}(\tau_1, \cdots, \tau_n) = 0, \quad w \in W_0, \alpha \in \mathbb{Z}^n,$$
where $p_{w,\alpha}(x)$ is the Laurent polynomial $\sum_\beta c_{w,\alpha,\beta} x^\beta$. Since $\tau_i(x^\gamma) = (q^{\gamma_i}) x^\gamma$, applying the expression to $x^\gamma$ we obtain
$$\sum x^{\alpha + w\gamma} p_{w,\alpha}(q^{\gamma_1}, \cdots, q^{\gamma_n}) = 0.$$

It follows that $p_{w,\alpha}(q^{\gamma_1}, \cdots, q^{\gamma_n}) = 0$ for all $\gamma$ in $\mathbb{Z}^n$ outside the *finite* union of hyperplanes determined by the conditions $\alpha + w\gamma = \alpha' + w'\gamma$ for $\alpha, \alpha', w, w'$ occuring in the last expression. But then $p_{w,\alpha}$ must be identically 0, and we conclude that all $c_{w,\alpha,\beta} = 0$, contrary to the assumption. $\square$

Let us define $\mathcal{R}_X := \mathbb{F}[X_1^{\pm 1}, \cdots, X_n^{\pm 1}]$; then the above proof shows:

3.3. COROLLARY.    *The natural maps from $H$, $H_0$, $\mathcal{R}_Y$ and $\mathcal{R}_X$ into $\mathcal{H}$ are injective.*

We shall identify the above algebras with their images in $\mathcal{H}$. Then we have:

3.4. COROLLARY.    *The multiplication maps from $\mathcal{R}_X \otimes H$ and $\mathcal{R}_X \otimes H_0 \otimes \mathcal{R}_Y$ into $\mathcal{H}$ are linear isomorphisms.*

We conclude this section by giving an intrinsic definition of the representation $\pi$. First, by the definition of $H$, it is clear that the map

(16) $$\chi : T_i \mapsto t_i^{1/2}, \quad i = 0, \cdots, n$$

extends to a one-dimensional character of $H$.

3.5. PROPOSITION. *The representation $\pi$ is isomorphic to $\mathrm{Ind}_H^{\mathcal{H}}(\chi)$.*

*Proof.* The induced representation is $\mathcal{H}/\mathcal{I}$ where $\mathcal{I}$ is the left ideal of $\mathcal{H}$ generated by $T_i - t_i^{1/2}, i = 0, \cdots, n$. On the other hand, $\pi \approx \mathcal{H}/\mathcal{J}$ where $\mathcal{J}$ is the annihilator of the cyclic vector $1 \in \mathcal{R}$. It remains to show that $\mathcal{I} = \mathcal{J}$.

First, since $s_i(1) = 1$ it follows from formula (7) that $\pi(T_i - t_i^{1/2})(1) = 0$, and so
$$\mathcal{I} \subseteq \mathcal{J}.$$

Next consider the left ideal $I$ in $H$ generated by $T_i - t_i^{1/2}$. Then we have $H = \mathbb{F} + I$; thus $\mathcal{R}_X \otimes H = \mathcal{R}_X \otimes \mathbb{F} + \mathcal{R}_X \otimes I$. Applying Corollary 3.4 we conclude that
$$\mathcal{H} = \mathcal{R}_X + \mathcal{I}.$$

This mean that the isomorphism $\mathcal{R}_X \to \mathcal{R}$ given by $X_i \mapsto x_i = \pi(X_i) \cdot 1$ can be factored as the following sequence of *surjective* maps
$$\mathcal{R}_X \to \mathcal{H}/\mathcal{I} \to \mathcal{H}/\mathcal{J} \to \mathcal{R}.$$

In particular, the middle map is an isomorphism, and so $\mathcal{I} = \mathcal{J}$. □

## 4. The involution

Let $\varepsilon$ denote the involution on $\mathbb{F}$ which sends $q, t, t_n, u_0$ to their inverses, and which maps
$$t_0 \mapsto u_n^{-1}.$$

We shall show that $\varepsilon$ extends to an involution on $\mathcal{H}$. First we prove

4.1. LEMMA. *Let $U_n \equiv X_1^{-1} T_0 Y_1^{-1} = X_1^{-1} T_1^{-1} \cdots T_n^{-1} \cdots T_1^{-1}$, then $U_n \sim u_n$.*

*Proof.* By (iv) we have $X_i^{-1} T_i^{-1} = T_i X_{i+1}^{-1}$, and applying this repeatedly, we get
$$U_n = (T_1 \cdots T_{n-1})(X_n^{-1} T_n^{-1})(T_{n-1}^{-1} \cdots T_1^{-1}).$$

Thus $U_n$ is conjugate to $X_n^{-1} T_n^{-1}$ and the lemma follows from relation (v). □

4.2. THEOREM. *The map $\varepsilon$ extends to an involution of $\mathcal{H}$ which maps $X_i$ to $Y_i$, sends $T_1, \cdots, T_n$ to their inverses, and maps*
$$T_0 \mapsto U_n^{-1}.$$

*Proof.* We first verify that the $\varepsilon$-transforms of (i)–(vi) hold in $\mathcal{H}$. For $i \neq 0$ the relation (i) becomes $T_i^{-1} \sim t_i^{-1}$, which is implied by $T_i \sim t_i$; while for $i = 0$ it becomes $U_n^{-1} \sim u_n^{-1}$, which follows from Lemma 4.1.

All the $\varepsilon$-transforms of the braid relations (ii) are immediate, except for
$$U_n^{-1}T_1^{-1}U_n^{-1}T_1^{-1} \stackrel{?}{=} T_1^{-1}U_n^{-1}T_1^{-1}U_n^{-1}.$$

We shall check this directly. Write $\Phi = T_2 \cdots T_n \cdots T_2$; then $X_1\Phi = \Phi X_1$, and we get

$$T_1^{-1}U_n^{-1}T_1^{-1}U_n^{-1} = \Phi T_1 X_1 \Phi T_1 X_1 = \Phi T_1 \Phi X_1 T_1 X_1 = \Phi T_1 \Phi X_1 X_2 T_1^{-1}$$
$$U_n^{-1}T_1^{-1}U_n^{-1}T_1^{-1} = T_1\Phi T_1 X_1 \Phi T_1 X_1 T_1^{-1} = T_1 \Phi T_1 \Phi X_1 T_1 X_1 T_1^{-1}$$
$$= T_1\Phi T_1 \Phi T_1^{-1} X_1 X_2 T_1^{-1}.$$

By multiplying both sides by $T_1 X_1^{-1} X_2^{-1} T_1$ on the right, it suffices to show
$$(T_2 \cdots T_n \cdots T_2)T_1(T_2 \cdots T_n \cdots T_2)T_1 \stackrel{?}{=} T_1(T_2 \cdots T_n \cdots T_2)T_1(T_2 \cdots T_n \cdots T_2).$$

We apply $T_2 T_1 T_2 = T_1 T_2 T_1$ to both sides and commute the resulting $T_1$'s as far to the extremes as possible. Using $T_1 T_2 T_1 = T_2 T_1 T_2$ once on each side and cancelling, we get

$$(T_3 \cdots T_n \cdots T_3)T_2(T_3 \cdots T_n \cdots T_3)T_2 \stackrel{?}{=} T_2(T_3 \cdots T_n \cdots T_3)T_2(T_3 \cdots T_n \cdots T_3).$$

Iterating, we reach the *true* relation $T_n T_{n-1} T_n T_{n-1} = T_{n-1} T_n T_{n-1} T_n$, which proves (ii).

The $\varepsilon$-transforms of (iii)–(iv) are easily checked. For (v) we get
$$\varepsilon(X_n^{-1}T_n^{-1}) = Y_n^{-1}T_n = (T_{n-1} \cdots T_1)T_0^{-1}(T_1^{-1} \cdots T_{n-1}^{-1})$$

which is conjugate to $T_0^{-1}$. Hence the desired relation follows from $T_0 \sim t_0$.

Finally, for (vi) we have
$$\varepsilon(U_0) = q^{1/2}U_n Y_1 = q^{1/2}(X_1^{-1}T_0 Y_1^{-1})Y_1 = q^{1/2}X_1^{-1}T_0 = U_0^{-1}.$$

Thus (vi) follows from its original counterpart.

It follows that $\varepsilon$ is a homomorphism, and it remains only to prove that $\varepsilon^2 = 1$. Since $\mathcal{H}$ is generated by $\{T_i, X_i, Y_i \mid i = 1, \cdots, n\}$, it suffices to show that $\varepsilon(Y_i) = X_i$ for all $i$. But
$$\varepsilon(Y_1) = T_1^{-1} \cdots T_n^{-1} \cdots T_1^{-1} U_n^{-1} = X_1.$$

Since $T_i X_i T_i = X_{i+1}$ and $T_i^{-1} Y_i T_i^{-1} = Y_{i+1}$, the result follows for all $i$ by induction. □

## 5. Intertwiners

In this section we introduce certain commutators in $\mathcal{H}$, and prove that they enjoy a crucial intertwining property with respect to the commutative family $\mathcal{R}_Y$.

*Definition.* We define operators $S_i$ in $\mathcal{H}$ as follows:

$$S_i := [T_i, Y_i], i = 1, \cdots, n; \quad S_0 := [Y_1, U_n]. \tag{17}$$

We also introduce the following notation, analogous to (7):

$$X^{v+k\delta} := q^{-k} X_1^{v_1} \cdots X_n^{v_n}; \quad Y^{v+k\delta} := q^k Y_1^{v_1} \cdots Y_n^{v_n}, \quad v \in \mathbb{Z}^n, \ k \in \mathbb{Z}. \tag{18}$$

5.1. THEOREM. *For all $\widetilde{v}$ in $\mathbb{Z}^n \times \mathbb{Z}\delta$,*

$$Y^{\widetilde{v}} S_i = S_i Y^{s_i(\widetilde{v})}; \quad i = 0, \cdots, n. \tag{19}$$

*Proof.* By Theorems 3.2 and 4.2, it is enough to prove the $\pi \circ \varepsilon$ transform of (19). Applying $\pi \circ \varepsilon$ to (18) and (17), we get:

$$\pi \circ \varepsilon(Y^{\widetilde{v}}) = x^{\widetilde{v}}; \quad \pi \circ \varepsilon(Y^{s_i(\widetilde{v})}) = x^{s_i(\widetilde{v})}$$

$$\pi \circ \varepsilon(S_0) = [x_1, \pi(T_0^{-1})]; \quad \pi \circ \varepsilon(S_i) = [\pi(T_i^{-1}), x_i], \ i = 1, \cdots, n.$$

Now, an easy calculation in $\text{End}(\mathcal{R})$, using formulas (13) gives

$$\pi \circ \varepsilon(S_i) = \phi_i(x) s_i; \quad \phi_i(x) = \begin{cases} t_0^{-1/2} x_1 (1 - cx_1^{-1})(1 - dx_1^{-1}) & i = 0 \\ t_n^{-1/2} x_n^{-1} (1 - ax_n)(1 - bx_n) & i = n \\ t_i^{-1/2} x_{i+1} (1 - t_i x_i x_{i+1}^{-1}) & i \neq 0, n. \end{cases} \tag{20}$$

Thus the $\pi \circ \varepsilon$ transform of (19) becomes the following assertion in $\text{End}(\mathcal{R})$:

$$x^{\widetilde{v}} \phi_i(x) s_i \stackrel{?}{=} \phi_i(x) s_i x^{s_i(\widetilde{v})}.$$

After cancelling the $\phi_i$, this follows from (8). $\square$

5.2. COROLLARY. *Let $a', b', c', d'$ be the $\varepsilon$ transforms of $a, b, c, d$, then*

$$S_i^2 = \begin{cases} u_n q^{-1} (1 - c' Y_1^{-1})(1 - d' Y_1^{-1})(1 - qc' Y_1)(1 - qd' Y_1) & i = 0 \\ t_n (1 - a' Y_n)(1 - b' Y_n)(1 - a' Y_n^{-1})(1 - b' Y_n^{-1}) & i = n \\ t_i Y_i Y_{i+1} (1 - t_i^{-1} Y_i Y_{i+1}^{-1})(1 - t_i^{-1} Y_i^{-1} Y_{i+1}) & i \neq 0, n. \end{cases} \tag{21}$$

*Proof.* By (20), we get $\pi \circ \varepsilon(S_i^2) = \phi_i s_i \phi_i s_i = \phi_i s_i (\phi_i) s_i^2 = \phi_i s_i (\phi_i)$, and the result follows from the explicit formula for $\phi_i$. $\square$

In the next section we will use the $S_i$'s as creation operators for the $E_\alpha$, starting with the constant function 1, which is an eigenfunction of $Y_i$, satisfying

$$\pi(Y_i)(1) = q^{\rho_i}(1); \quad i = 1, \cdots, n. \tag{22}$$

This is an immediate consequence of the equation $\pi(T_i)(1) = t_i^{1/2}(1)$ and the definitions of $Y_i$ and $\rho$ in (11) and (15). To describe the other eigenvalues, we proceed as follows:

*Definition.* For $\alpha$ in $\mathbb{Z}^n$, we define

$w_\alpha :=$ the shortest element in $W_0$ such that $w_\alpha^{-1} \cdot \alpha$ is a partition;

$\overline{\alpha} := \alpha + w_\alpha \cdot \rho$ where $\rho$ is as in (15);

$\mathcal{R}_\alpha :=$ the space of all $f \in \mathcal{R}$ satisfying $Y_i f = q^{\alpha_i + (w_\alpha \cdot \rho)_i} f$ for all $i$.

Alternatively, $w_\alpha$ in $W_0 = (\pm 1)^n S_n$ can be described as $w_\alpha := \sigma_\alpha \pi_\alpha$, where $\sigma_\alpha \in (\pm 1)^n$ is simply $(\mathrm{sgn}(\alpha_1), \cdots, \mathrm{sgn}(\alpha_n))$ with $\mathrm{sgn}(0)$ defined to be 1; and $\pi_\alpha$ is the permutation in $S_n$ defined as follows: order the indices *first* by decreasing $|\alpha_i|$, *then* for fixed $|\alpha_i|$ from left to right for $\alpha_i \geq 0$, and *finally* from right to left for $\alpha_i < 0$.

For example if $\alpha = (-2, 2, 1, -1, 0, 1, -1)$, then $\sigma_\alpha = (-1, 1, 1, -1, 1, 1, -1)$ and $\pi_\alpha$ is the permutation $(2, 1, 3, 6, 7, 4, 5)$.

5.3. THEOREM. *If $s_i \cdot \alpha \neq \alpha$, then $\pi(S_i)$ is a linear isomorphism from $\mathcal{R}_\alpha$ to $\mathcal{R}_{s_i \cdot \alpha}$.*

*Proof.* Let $f \in \mathcal{R}_\alpha$. Then by (5) and (18) it follows that for all $\widetilde{v}$ in $\mathbb{Z}^n \times \mathbb{Z}\delta$,
$$Y^{\widetilde{v}}(f) = q^{\langle \widetilde{v}, \alpha + w_\alpha \cdot \rho \rangle} f.$$

Let us write $\overline{\alpha} = \alpha + w_\alpha \cdot \rho$. Then from Theorem 5.1 and (3) we get
$$\pi(Y^{\widetilde{v}})\pi(S_i) f = \pi(S_i)\pi(Y^{s_i(\widetilde{v})}) f = q^{\langle s_i(\widetilde{v}), \overline{\alpha} \rangle} \pi(S_i) f = q^{\langle \widetilde{v}, s_i \cdot \overline{\alpha} \rangle} \pi(S_i) f.$$

Thus to prove that $\pi(S_i) f \in \mathcal{R}_{s_i \cdot \alpha}$, it suffices to show that $s_i \cdot \alpha \neq \alpha$ implies

(23) $$s_i \cdot \overline{\alpha} = \overline{s_i \cdot \alpha}; \quad i = 0, \cdots, n.$$

For $i = 0$, write $\beta = s_0 \cdot \alpha = (-\alpha_1 - 1, \alpha_2, \cdots, \alpha_n)$. Then we claim that the permutations $\pi_\beta$ and $\pi_\alpha$ are the *same*. Indeed if $\alpha_1$ is positive then, in the ordering corresponding to $\pi_\alpha$, the index 1 is the first among the indices $j$ with $|\alpha_j| = \alpha_1$; while in the ordering for $\pi_\beta$, 1 is the last index in the higher group of indices satisfying $|\beta_j| = \alpha_1 + 1$. Thus the relative position of 1 with respect to other indices stays the same, as do those of other indices with respect to each other. A similar argument works if $\alpha_1$ is negative, and taking into account the sign change we conclude that:

$$\overline{\beta}_1 = \beta_1 + (w_\beta \cdot \rho)_1 = -\alpha - 1 - (w_\alpha \cdot \rho)_1 = -\overline{\alpha}_1 - 1; \quad \overline{\beta}_i = \overline{\alpha}_i,\ i \geq 2,$$

which is precisely the content of (23) for $i = 0$.

The argument for $i \geq 1$ is similar and simpler. We observe that if $\beta := s_i \cdot \alpha$ is different from $\alpha$, then $w_\beta$ equals $s_i w_\alpha$. Since $s_i$ acts linearly, we get
$$s_i \cdot \overline{\alpha} = s_i \cdot \alpha + s_i \cdot (w_\alpha \cdot \rho) = \beta + w_\beta \cdot \rho = \overline{\beta} = \overline{s_i \cdot \alpha}.$$

Thus we conclude that $\pi(S_i)$ maps $\mathcal{R}_\alpha$ into $\mathcal{R}_{s_i\cdot\alpha}$ for all $i \geq 0$. But, by Corollary 5.2, $S_i^2$ is in $\mathcal{R}_Y$; hence $\pi(S_i^2)$ acts by a scalar $c_i$ on $\mathcal{R}_\alpha$, which can be readily computed by substituting $Y_i = \overline{\alpha}_i$ in (21). In particular, we see that if $s_i \cdot \alpha \neq \alpha$ then $c_i$ is not zero. Thus $\pi(S_i)$ is a linear isomorphism from $\mathcal{R}_\alpha$ to $\mathcal{R}_{s_i\cdot\alpha}$, with inverse $c_i^{-1}\pi(S_i)$. □

## 6. Nonsymmetric Koornwinder polynomials

In this section we will define the nonsymmetric Koornwinder polynomials. The crucial result is:

6.1. THEOREM. *The spaces $\mathcal{R}_\alpha$ are all one-dimensional.*

*Proof.* We first prove that the spaces $\mathcal{R}_\alpha$ are nonzero. For $\alpha = 0$, $w_\alpha$ is the identity in $W_0$ and so $\overline{\alpha} = \rho$. Thus by (22), the constant functions belong to $\mathcal{R}_\alpha$ for $\alpha = 0$. For other $\alpha \in \mathbb{Z}^n$ we use Theorem 5.3 together with the fact that the affine action of $W$ on $\mathbb{Z}^n$ is transitive.

Now let $f = \sum_{c_\beta \in \mathbb{F}} c_\beta x^\beta$ be a nonzero function in $\mathcal{R}_\alpha$. Then $f$ satisfies
$$\pi(Y_i)f := q^{(\alpha+w_\alpha\cdot\rho)_i}f; \quad i = 1, \cdots, n.$$
As in the proof of Theorem 3.2, we set $t, t_0, t_n, u_0, u_n$ equal to 1 in the expression. Then $\rho$ specializes to the zero vector, and $Y_i$ specializes to $\pi(\tau_i) = T_{q,x_i}$. Clearing denominators and eliminating common factors, we may also assume that the $c_\beta$ have finite specializations, not all zero. Letting $g \neq 0$ denote the specialization of $f$, we get
$$T_{q,x_i}g = q^{\alpha_i}g$$
which means that $g$ is a nonzero multiple of $x^\alpha$.

In particular, the coefficient $c_\alpha$ has a nonzero specialization and so must be nonzero. The result follows, since if there were two linearly independent functions in $\mathcal{R}_\alpha$, we could construct a nonzero $f$ with $c_\alpha = 0$. □

The proof of the theorem shows that a function $f$ in $\mathcal{R}_\alpha$ is uniquely determined by the knowledge of the coefficient of $x^\alpha$ in $f$.

*Definition.* The nonsymmetric Koornwinder polynomial $E_\alpha$ is the unique polynomial in $\mathcal{R}_\alpha$ in which the coefficient of $x^\alpha$ is 1.

6.2. THEOREM. *The polynomials $E_\alpha$ form a basis for $\mathcal{R}$ over $\mathbb{F}$.*

*Proof.* Let us consider the *degree* filtration $\mathcal{R}_{(0)} \subseteq \mathcal{R}_{(1)} \subseteq \cdots \subseteq \mathcal{R}$, where $\mathcal{R}_{(k)}$ is spanned by all monomials $x^\alpha$ with $|\alpha| := |\alpha_1| + \cdots + |\alpha_n| \leq k$. We claim that for $|\alpha| \leq k$,
$$\mathcal{R}_\alpha \subseteq \mathcal{R}_{(k)}.$$

For $\alpha = 0$, we use Theorem 6.1 to see that $\mathcal{R}_\alpha$ consists precisely of the constants, which lie in $\mathcal{R}_{(0)}$. For other $\alpha$, we observe that by (13) the filtration is invariant under $T_i$ and $Y_i$, while $U_n = X_1^{-1} T_0 Y_1^{-1}$ raises degree by at most one. Thus $S_1, \cdots, S_n$ preserve the filtration while $S_0$ raises degree by at most one. Now any $\alpha$ can be obtained from 0 by applying a sequence of $s_i$'s in which $s_0$ occurs exactly $|\alpha|$ times. Applying the corresponding $S_i$'s to $\mathcal{R}_0$, we deduce the claim.

It follows that the set $\{E_\alpha : |\alpha| \leq k\}$ is contained in $\mathcal{R}_{(k)}$ and has the same cardinality as the monomial basis. Therefore it suffices to prove that the $E_\alpha$ are linearly independent. For this we choose a polynomial $f$ in $\mathcal{R}$ which takes distinct values on the finite set $\{q^{\overline{\alpha}} : |\alpha| \leq k\}$. Then the $E_\alpha$ belong to distinct eigenspaces under the operator $\pi(f(Y_1, \cdots, Y_n))$ and hence are linearly independent. $\square$

6.3. COROLLARY. *The representation $\pi$ is irreducible.*

*Proof.* Let $\mathcal{V}$ be a $\pi(\mathcal{H})$-invariant subspace of $\mathcal{R}$, and suppose $f = \sum c_\alpha E_\alpha$ belongs to $\mathcal{V}$, with some $c_\beta \neq 0$. Choose a function $g$ in $\mathcal{R}$ such that $g(q^{\overline{\beta}}) = 1/c_\beta$, and $g(q^{\overline{\alpha}}) = 0$ for all other $\alpha$ for which $c_\alpha \neq 0$. Then applying $\pi(g(Y_1, \cdots, Y_n))$ to $f$ we conclude that $E_\beta$ belongs to $\mathcal{V}$. Now applying the $\pi(S_i)$'s we conclude that every $E_\alpha$ belongs to $\mathcal{V}$. $\square$

Next we consider the restriction of $\pi$ to $H$.

*Definition.* For each partition $\lambda$ we write $\mathcal{R}^\lambda$ for the subspace of $\mathcal{R}$ spanned by the $\{E_\alpha : \alpha \in W_0 \cdot \lambda\}$.

6.4. COROLLARY. *The $\mathcal{R}^\lambda$ are irreducible $\pi(H)$-modules, and $\mathcal{R}$ is their direct sum.*

*Proof.* For the irreducibility, we repeat the previous argument without involving $S_0$. The second assertion follows from Theorem 6.2. $\square$

Finally we discuss the connection with the symmetric Koornwinder polynomials $P_\lambda$.

6.5. COROLLARY. *The symmetric Koornwinder polynomial $P_\lambda$ can be characterized as the unique $W_0$-invariant polynomial in $\mathcal{R}^\lambda$ which has the coefficient of $x^\lambda$ equal to 1.*

*Proof.* We note that if $\alpha \in W_0 \cdot \lambda$, then $\overline{\alpha} = \alpha + w_\alpha \cdot \rho = w_\alpha \cdot (\lambda + \rho)$. In particular, if $f$ is in $\mathcal{S}$, then $f(\overline{\alpha}) = f(\lambda + \rho)$, and so $\mathcal{R}^\lambda$ is precisely the $f(\lambda + \rho)$-eigenspace of $\pi(f(Y_1, \cdots, Y_n))$ for $f \in \mathcal{S}$. The result now follows from the characterization (14). $\square$

*Definition.* Define $C \in H_0$ by $C := \left(\sum_{w \in W_0} \chi(T_w)^2\right)^{-1} \sum_{w \in W_0} \chi(T_w) T_w$.

6.6. COROLLARY.   $\pi(C)$ *is a projection from* $\mathcal{R}^\lambda$ *to* $\mathbb{F} P_\lambda$.

*Proof.* First of all, an easy calculation as in Lemma 2.5 of [S] shows that $T_i C = t_i^{1/2} C$ for $i = 1, \cdots, n$; hence $\pi(T_i)\pi(C)f = t_i^{1/2}\pi(C)f$ for all $f \in \mathcal{R}$. By (13), this implies that $\pi(C)f$ is $W_0$-invariant, and so must be a multiple of $P_\lambda$.

Moreover, for $f \in \mathcal{S}$, $\pi(T_w)f = \chi(T_w)f$; hence $\pi(C)$ acts by the identity on $\mathcal{S}$. $\square$

## 7. Duality

Let $^\dagger$ denote the involution on $\mathbb{F}$ which sends $q, t, t_0, t_n, u_0, u_n$ to their inverses.

7.1. PROPOSITION.   *The map* $^\dagger$ *extends to an anti-involution on* $\mathcal{H}$ *satisfying*
$$T_i^\dagger = T_i^{-1},\ X_i^\dagger = X_i^{-1},\ Y_i^\dagger = Y_i^{-1}.$$

*Proof.* For the proof we merely observe that each defining relation of $\mathcal{H}$ is $^\dagger$-invariant. $\square$

*Definition.* We define the *duality* anti-involution $^*$ on $\mathcal{H}$ by
$$h^* = \varepsilon(h^\dagger) = \varepsilon(h)^\dagger, \quad h \in \mathcal{H}.$$

On $\mathbb{F}$, $^*$ simply switches $t_0$ and $u_n$; while on the generators,
$$T_i^* = T_i,\ X_i^* = Y_i^{-1},\ Y_i^* = X_i^{-1},\ i = 1\cdots, n; \quad T_0^* = U_n.$$

We also extend $^*$ from $\mathbb{F}$ to an involution on $\mathcal{R}$ by defining $x_i^* = x_i^{-1}$ for all $i$. Observe that if $f$ is in $\mathcal{S}$, then $f$ is invariant under $x_i \mapsto x_i^{-1}$, and so $f^*$ is obtained just by switching $t_0$ and $u_n$ in the coefficients of $f$.

Next, we define $\rho^*$ by the requirement that $q^{\rho^*} = (q^\rho)^*$. Explicitly,
$$q^{\rho_i^*} = (q^{\rho_i})^* = ((t_0 t_n)^{1/2} t^{n-i})^* = (u_n t_n)^{1/2} t^{n-i}.$$

The duality conjecture of Macdonald can be stated as follows:

7.2. CONJECTURE.   *For any two partitions* $\lambda$ *and* $\mu$ *we see that*
$$\frac{P_\lambda(q^{\mu+\rho^*})}{P_\lambda(q^{\rho^*})} = \frac{P_\mu^*(q^{\lambda+\rho})}{P_\mu^*(q^\rho)}.$$

This is seen to be equivalent to the formulation in (4.4) of [D], after the easy verification that our definition of duality $(t_0 \leftrightarrow u_n)$, is the same as that in (4.1) of [D].

To establish Conjecture 7.2 and its analog for $E_\alpha$, we introduce the following:

*Definition.* Let $S$ be the map from $\mathcal{H}$ to $\mathbb{F}$ defined by
$$S(h) := F_h(q^{-\rho^*}); \quad \text{where } F_h = \pi(h)(1) \in \mathcal{R}.$$

7.3. THEOREM. *We have $S(h^*) = S(h)^*$ for all $h \in \mathcal{H}$.*

*Proof.* By linearity and Corollary 3.4 it is enough to prove this for $h$ of the form $X^\alpha T_w Y^\beta$, with $\alpha, \beta \in \mathbb{Z}^n$, $w \in W_0$. Then by (16) and (22),

(24) $\quad F_h = q^{\langle \beta, \rho \rangle} \chi(T_w) x^\alpha; \quad \text{and} \quad S(h) = q^{\langle \beta, \rho \rangle} \chi(T_w) q^{-\langle \alpha, \rho^* \rangle}.$

On the other hand $h^* = X^{-\beta} T_w^* Y^{-\alpha}$, and so

(25) $\quad F_{h^*} = q^{-\langle \alpha, \rho \rangle} \chi(T_w^*) x^{-\beta}; \quad \text{and} \quad S(h^*) = q^{-\langle \alpha, \rho \rangle} \chi(T_w^*) q^{\langle \beta, \rho^* \rangle}.$

For $w \in W_0$, $\chi(T_w)$ only involves $t$ and $t_n$; and $T_w^*$ is simply obtained from $T_w$ by reversing its product expansion (10) in terms of $T_1, \cdots, T_n$. Thus
$$\chi(T_w)^* = \chi(T_w) = \chi(T_w^*); \quad w \in W_0.$$

Since $*$ interchanges $\rho$ and $\rho^*$, the result now follows by comparing (24) and (25). □

We now define scalars $\mathcal{E}_{\alpha\beta}, \mathcal{P}_{\lambda\mu}$ in $\mathbb{F}$ by
$$\mathcal{E}_{\alpha\beta} := E_\alpha^*(q^{\overline{\beta}}) E_\beta(q^{-\rho^*}); \quad \mathcal{P}_{\lambda\mu} := P_\lambda^*(q^{\mu+\rho}) P_\mu(q^{-\rho^*}).$$

7.4. THEOREM. *We have $\mathcal{E}_{\alpha\beta}^* = \mathcal{E}_{\beta\alpha}$ and $\mathcal{P}_{\lambda\mu}^* = \mathcal{P}_{\mu\lambda}$.*

*Proof.* For the first assertion we consider $h := E_\alpha^*(Y) E_\beta(X)$. Then by the definition of $*$, we get $h^* = E_\beta^*(Y) E_\alpha(X)$. Now,
$$F_h := \pi(E_\alpha^*(Y)) E_\beta(x) = E_\alpha^*(q^{\overline{\beta}}) E_\alpha(x),$$
and so $S(h) = \mathcal{E}_{\alpha\beta}$. Similarly $S(h^*) = \mathcal{E}_{\beta\alpha}$, and the result follows from Theorem 7.3.

The second assertion is proved similarly by considering $h := P_\lambda^*(Y) P_\mu(X)$. □

7.5. COROLLARY. *Conjecture 7.2 is true.*

*Proof.* Since $P_\lambda$ is invariant under $x_i \mapsto x_i^{-1}$ we get
$$\mathcal{P}_{\mu\lambda} := P_\mu^*(q^{\lambda+\rho}) P_\lambda(q^{-\rho^*}) = P_\mu^*(q^{\lambda+\rho}) P_\lambda(q^{\rho^*}).$$

On the other hand,
$$\mathcal{P}_{\lambda\mu}^* := \left( P_\lambda^*(q^{\mu+\rho}) P_\mu(q^{-\rho^*}) \right)^* = P_\lambda(q^{-\mu-\rho^*}) P_\mu^*(q^\rho) = P_\lambda(q^{\mu+\rho^*}) P_\mu^*(q^\rho).$$

Thus the result follows from Theorem 7.4. □

Rutgers University, New Brunswick, NJ
*E-mail address*: sahi@math.rutgers.edu